\documentclass[a4paper,12pt]{article}

\pdfoutput=1

\usepackage[paper=letterpaper,margin=1in]{geometry}

\usepackage{amsmath,amssymb,amsfonts,epsfig,cite,setspace,bigstrut,longtable,array,breqn,color}
\usepackage{hyperref,caption,setspace}
\usepackage{amsthm}
\usepackage{pdfpages,lipsum}
\usepackage{todonotes}
\usepackage{cancel}

\newenvironment{red}{\relax\color{red}}{\relax}
\newenvironment{blue}{\relax\color{blue}}{\hspace*{.5ex}\relax}
\newcommand{\ber}{\begin{red}}
\newcommand{\er}{\end{red}}
\newcommand{\beb}{\begin{blue}}
\newcommand{\eb}{\end{blue}}

\newcommand{\sage}{{\sc SageMath }} 
\newcommand{\mth}{{\sc Mathematica }}
\begin{document}


\vskip 0.25in

\newcommand{\nn}{\nonumber}
\newcommand{\tr}{\mathop{\rm Tr}}
\newcommand{\comment}[1]{}
\newcommand{\cM}{{\cal M}}
\newcommand{\cW}{{\cal W}}
\newcommand{\cN}{{\cal N}}
\newcommand{\cH}{{\cal H}}
\newcommand{\cK}{{\cal K}}
\newcommand{\cZ}{{\cal Z}}
\newcommand{\cO}{{\cal O}}
\newcommand{\cA}{{\cal A}}
\newcommand{\cB}{{\cal B}}
\newcommand{\cC}{{\cal C}}
\newcommand{\cD}{{\cal D}}
\newcommand{\cT}{{\cal T}}
\newcommand{\cV}{{\cal V}}
\newcommand{\cE}{{\cal E}}
\newcommand{\cF}{{\cal F}}
\newcommand{\cX}{{\cal X}}
\newcommand{\IA}{\mathbb{A}}
\newcommand{\IP}{\mathbb{P}}
\newcommand{\IQ}{\mathbb{Q}}
\newcommand{\IH}{\mathbb{H}}
\newcommand{\IR}{\mathbb{R}}
\newcommand{\IC}{\mathbb{C}}
\newcommand{\IF}{\mathbb{F}}
\newcommand{\IV}{\mathbb{V}}
\newcommand{\II}{\mathbb{I}}
\newcommand{\IZ}{\mathbb{Z}}
\newcommand{\re}{{\rm~Re}}
\newcommand{\im}{{\rm~Im}}

\let\oldthebibliography=\thebibliography
\let\endoldthebibliography=\endthebibliography
\renewenvironment{thebibliography}[1]{%
\begin{oldthebibliography}{#1}%
\setlength{\parskip}{0ex}%
\setlength{\itemsep}{0ex}%
}%
{%
\end{oldthebibliography}%
}

\newtheorem{theorem}{\bf THEOREM}
\def\thetheorem{\thesection.\arabic{theorem}}
\newtheorem{proposition}{\bf PROPOSITION}
\def\thetheorem{\thesection.\arabic{proposition}}
\newtheorem{observation}{\bf OBSERVATION}
\def\thetheorem{\thesection.\arabic{observation}}
\newtheorem{conjecture}{\bf CONJECTURE} 
\def\thetheorem{\thesection.\arabic{CONJECTURE}}

\theoremstyle{definition}
\newtheorem{definition}{\bf DEFINITION} 
\def\thetheorem{\thesection.\arabic{DEFINITION}}
\newtheorem{example}{\bf EXAMPLE} 
\def\thetheorem{\thesection.\arabic{EXAMPLE}}
\newtheorem{remark}{\bf REMARK} 
\def\thetheorem{\thesection.\arabic{REMARK}}

\def\theequation{\thesection.\arabic{equation}}
\newcommand{\setall}{\setcounter{equation}{0}
        \setcounter{theorem}{0}}
\newcommand{\setequation}{\setcounter{equation}{0}}
\renewcommand{\thefootnote}{\fnsymbol{footnote}}

\newcommand{\seteq}{\mathbin{:=}}
\newcommand{\GL}{\operatorname{GL}}
\newcommand{\Sp}{\operatorname{Sp}}
\newcommand{\USp}{\operatorname{USp}}
\newcommand{\GSp}{\operatorname{GSp}}
\newcommand{\U}{\operatorname{U}}
\newcommand{\SU}{\operatorname{SU}}
\newcommand{\SO}{\operatorname{SO}}
\newcommand{\End}{\operatorname{End}}

\begin{titlepage}

~\\
\vskip 1cm

\begin{center}
{\Large \bf Machine-Learning Number Fields}
\end{center}
\medskip

\renewcommand{\arraystretch}{0.5} 

\vspace{.4cm}
\centerline{
{\large Yang-Hui He, Kyu-Hwan Lee, Thomas Oliver}
}
\vspace*{3.0ex}

\vspace{10mm}

\begin{abstract}
We show that standard machine-learning algorithms may be trained to predict certain invariants of algebraic number fields to high accuracy. A random-forest classifier that is trained on finitely many Dedekind zeta coefficients is able to distinguish between real quadratic fields with class number $1$ and $2$, to $0.96$ precision. Furthermore, the classifier is able to extrapolate to fields with discriminant outside the range of the training data. When trained on the coefficients of defining polynomials for Galois extensions of degrees $2$, $6$, and $8$, a logistic regression classifier can distinguish between Galois groups and predict the ranks of unit groups with precision $>0.97$.
\end{abstract}

\end{titlepage}

\begin{spacing}{1}
\tableofcontents
\end{spacing}

\section{Introduction \& Summary}

Algebraic number fields are characterized by various invariants. One such invariant is the class number, which encodes how far the ring of integers in the number field is from being a unique factorisation domain. To this day, there remain central open questions regarding class numbers of algebraic number fields. For example, whilst there is a well-known list of imaginary quadratic number fields with class number $1$, it is not known whether or not there are infinitely many real quadratic fields with class number $1$. Actually, Gauss conjectured in his famous {\em Disquisitiones Arithmeticae} of 1801 that there are infinitely many, and the Cohen--Lenstra heuristics predict that around 75\% of real quadratic fields would have class number $1$ \cite{CL}. In this paper, we show that a machine-learning algorithm may be trained to predict certain invariants, including the class number of real quadratic fields.

For a broad introduction to machine-learning, see \cite{ML,Hastie}. The machine-learning of mathematical structures is a relatively recent enterprise. Interesting early neural-network experiments exploring the non-trivial zeros of the Riemann zeta function were documented in \cite{shanker} (a more recent work is \cite{KV}). Building on work in superstring theory, more precisely the computation of topological invariants for Calabi--Yau compactifications \cite{He:2017set,Krefl:2017yox,Ruehle:2017mzq,Carifio:2017bov} (q.v., \cite{He:2018jtw} for a summary), a programme developing the applications of machine-learning to abstract mathematics was proposed in \cite{He:2017aed,He:2017set}. Since then, machine-learning has been applied to various branches within the discipline with the intention of pattern-recognition and conjecture-raising. To name a few: representation theory \cite{He:2019nzx}, graph theory \cite{He:2020fdg}, metric geometry \cite{Ashmore:2019wzb}, dessins d'enfants \cite{He:2020eva}, and quiver mutations \cite{Bao:2020nbi}.  

Recently, the present authors demonstrated that techniques from machine-learning could be used to resolve a classification problem in arithmetic geometry \cite{HLO}. To be precise, we showed that a Bayesian classifier can distinguish between Sato--Tate groups given a small number of Euler factors for the $L$-function with over $99\%$ accuracy. Given the efficient nature of the machine-learning approach, [Loc.~cit.] suggested a machine can be trained to learn the Sato--Tate distributions and may be able to classify curves much more efficiently than the methods available in the literature.

This paper is a continuation of our observation that machine-learning can be used in number theory. In particular, we will apply logistic regression and random forest classifiers---these are reviewed in \cite[Sections~4.4~\&~15]{Hastie}. Our experiments are concerned with predicting the following invariants: degree, signature, Galois group, and class number. We utilise three training data sets associated to algebraic number fields: (1) coefficients of their defining polynomials, (2) finitely many coefficients of their Dedekind zeta functions, and (3) binary vectors encoding finitely many completely-split rational primes. Each training dataset has its own strengths and weaknesses. We review the utility of datasets (1), (2) and (3) below. 

\begin{enumerate}
\item Using both defining polynomial and zeta coefficient training, we observe high-accuracy predictions for number field signatures and Galois groups. The signature of a number field determines the rank of its unit group, which is equal to the vanishing order of the associated Dedekind zeta function at $s=0$. Elsewhere it has been demonstrated that a machine cannot be efficiently trained to predict the ranks of elliptic curves from their minimal Weierstrass equation \cite{Alessandretti:2019jbs}, which is in contrast to our observations for number fields. 
\item The Dedekind zeta function of an algebraic number field has a simple pole at $s=1$ with residue given by the analytic class number formula. 
We train a random forest classifier through 1000 zeta coefficients of real quadratic fields with class number $1$ or $2$ with discriminant less than one million which are available at \cite[Number~Fields]{lmfdb}, and the resulting classifier can distinguish between class numbers $1$ and $2$ with accuracy $0.96$. When we apply the same classifier to real quadratic fields with discriminants between $1$ million and $3$ million, we find that it distinguish between class numbers $1$ and $2$ with accuracy $0.92$.
\item It is well-known that the set of split primes uniquely characterizes a Galois extension over $\mathbb Q$, cf. \cite[VII, \S 13]{Neu}. Motivated by this, we train classifiers using binary data recording split primes and apply the classifiers to various invariants of number fields. However, the classifiers perform poorly except for detecting degrees of the extensions.  
\end{enumerate}

An outline of the contents of this paper is as follows. In Section~\ref{s:nomen} we recall basic terminology and establish the notation used in the sequel. In Section~\ref{sec:TrainingData} we define our three forms of training data, and explain the experimental set-up. In Section~\ref{sec:signature} it is shown that, when trained on zeta coefficients, a random forest classifier is able to distinguish between extension degrees and signatures. Furthermore, we apply logistic regression to the defining polynomial dataset. In Section~\ref{sec:GalGroups} we outline our experiments with Galois groups of order $8$. In this case, zeta coefficients and defining polynomial coefficients perform equally well. In Section~\ref{sec:classnumbers} it is observed that, when trained on zeta coefficients, a random forest classifier is able to distinguish between real quadratic fields of class number $1$ and class number $2$. The classifier is trained using quadratic fields with discriminant less than one million, but is able to extrapolate to ranges far beyond the training data.

\subsection*{Acknowledgements}
YHH is indebted to STFC UK, for grant ST/J00037X/1, KHL is partially supported by a grant from the Simons Foundation (\#712100), and TO acknowledges support from the EPSRC through research grant EP/S032460/1.

\section{Nomenclature}\label{s:nomen}

We will use the following notation throughout:

\begin{description}

\item[Algebraic number field] denoted by $F$. We will assume that the extension $F/\mathbb{Q}$ is Galois;

\item[Extension degree] of $F$ over $\IQ$ is denoted $[F : \IQ ]$;

\item[Signature] of $F$ is the pair $(r_1,r_2)$, in which $r_1$ (resp.~$r_2$) denotes the number of real embeddings (resp.~conjugate pairs of complex embeddings) of $F$.  If $(r_1,r_2)$ is the signature of $F$, then $[F:\mathbb{Q}]=r_1+2r_2$. If $r_2=0$ (resp.~$r_1=0$) then we refer to $F$ as totally real (resp.~imaginary);

\item[Ring of integers] denoted by $\mathcal{O}_F$;

\item[Rank] of the unit group $\mathcal O^\times_F$ is equal to $r:=r_1+r_2-1$ by Dirichlet's unit theorem;

\item[Discriminant] of $F$ denoted by $\Delta_F$; it is known that $\mathrm{sgn}(\Delta_F)=(-1)^{r_2}$;

\item[Ramification] A rational prime $p$ {\bf ramifies} in $F$ if and only if $p$ divides $\Delta_F$; an unramified prime $p$ {\bf splits completely} in $F$ if $p \mathcal O_F$ is a product of $[F:\mathbb Q]$-many distinct prime ideals in $\mathcal O_F$, and $p$ is {\bf inert} in $F$ if $p \mathcal O_F$ is itself a prime ideal;

\item[Class number] of $F$ denoted by $h_F$. That is, the size of the ideal class group (the quotient group of the fractional ideals by the principal ideals);

\item[Norm] of an ideal $I$ in $\mathcal{O}_F$ is denoted by $N(I)$;

\item[Prime ideal] denoted by $\mathfrak{p}$. A prime ideal ideal in $\mathcal{O}_F$ lies above a rational prime $p$ if $\mathfrak{p}$ divides the ideal generated by $p$; we denote this situation by $\mathfrak{p}|p$;

\item[Quadratic number field] has the form $\mathbb{Q}(\sqrt{d})$ with $d$ a square-free integer.
If $d<0$ (resp.~$d>0$) then we call the field imaginary quadratic (resp.~real quadratic). 
The discriminant of $F=\mathbb{Q}(\sqrt{d})$ is $d$ (resp. $4d$) if $d\equiv1\text{ mod }4$ (resp.~$d\equiv2,3\text{ mod }4$). In particular, a real quadratic number field has positive discriminant;

\item[Galois group] associated to the Galois extension $F/\mathbb{Q}$ is denoted by $\mathrm{Gal}(F/\mathbb{Q})$;

\item[Cyclic group] of order $n$ denoted by $C_n$;

\item[Dihedral group] of order $2n$ denoted by $D_n$.

\end{description}

\section{Establishing the Datasets}\label{sec:TrainingData}

In this section we explain our training datasets, and outline the basic experimental strategy.

\subsection{Defining polynomials}

Recall from Section~\ref{s:nomen} that we assume the extension $F/\mathbb{Q}$ to be Galois. A defining polynomial for $F$ is an irreducible polynomial $P(x)\in\mathbb{Q}[x]$ such that $F=\mathbb{Q}(\alpha)$ for a root $\alpha$ of $P(x)$. We choose $P(x)$ as in \cite[Normalization~of~defining~polynomials~for~number~fields]{lmfdb}. In particular, $P(x)$ is monic with integer coefficients, and, if $\alpha_1,\dots,\alpha_n$ are the complex roots of $P(x)$, then the sum $\sum_{i=1}^n|\alpha_i|^2$ is minimized. We write:
\begin{equation}\label{eq.minpoly}
P(x)=x^n+c_{n-1}x^{n-1}+\cdots+c_1x+c_0, \ \ c_i\in\mathbb{Z}, \ \ n=[F:\mathbb{Q}].
\end{equation}
Using the coefficients of $P(x)$, we define the vector:
\begin{equation}\label{eq.polyvector}
v_P(F)=(c_0,\dots,c_{n-1})\in\mathbb{Z}^n.
\end{equation}
Let $\mathcal{F}$ denote a finite set of number fields, and, for all $F\in\mathcal{F}$, let $c(F)$ be an invariant of interest. For example, $\mathcal{F}$ could be the set of all real quadratic fields with discriminant less than one million and, for $F\in\mathcal{F}$, the invariant $c(F)$ could be the class number of $F$. We introduce the following labeled dataset: 
\begin{equation}\label{eq.polydata}
\mathcal{D}_P=\{v_P(F)\rightarrow c(F):F\in\mathcal{F}\}.
\end{equation}
\begin{example}\label{ex.6206}
In Section~\ref{ss:abnonab}, we will take $\mathcal F$ to contain certain degree 8 number fields with Galois group isomorphic to either $C_8$ or $D_4$. For $F\in\mathcal{F}$ we will let $c(F)$ be $0$ (resp. $1$) corresponding to $\mathrm{Gal}(F/\mathbb{Q})\cong C_8$ (resp. $\mathrm{Gal}(F/\mathbb{Q})\cong D_4$). A large database of such fields can be downloaded from \cite[Number~Fields]{lmfdb}, including around $6200$ such that $c(F)=0$. The set $\mathcal{F}$ consists of these fields, and a random sample of around $6200$ (out of around $28000$) fields such that $c(F)=1$.  An instance of $v_P(F)$ such that $c(F)=0$ is 
\[ 
(4096, -512, 320, 136, -46, 17, 5, -1).
\]
The $v_P(F)$ with the largest $c_0$ such that $c(F)=0$ is
\begin{align*} 
&(153220409851123184, 812631532526484, 13364221512257,\\ & \phantom{LLLLLLLLL} -78983668469, -234643970, -7256689, 11478, -1).  
\end{align*}
\end{example}

\subsection{Dedekind zeta functions}

The Dedekind zeta function of a number field $F$ is given by the following formulas:
\[
\zeta_F(s)=\prod_{\mathfrak{p}}\left(1-N(\mathfrak{p})^{-s}\right)^{-1}=\sum_{I\leq\mathcal{O}_F}N(I)^{-s}=\sum_{n=1}^{\infty}a_nn^{-s},
\]
where $\mathfrak{p}$ varies over prime ideals in $\mathcal{O}_F$, $I$ varies over the non-zero ideals in $\mathcal{O}_F$, and, for a positive integer $n$, 
\begin{equation}\label{eq.an}
a_n=\#\{N(I)=n:I\leq\mathcal{O}_F\},\ \ n\in\mathbb{Z}_{\geq1}.
\end{equation}
Since we assume that $F$ is Galois over $\mathbb Q$, the zeta function $\zeta_F(s)$ uniquely determines $F$. However, we caution that in general a number field is not determined by its Dedekind zeta function\footnote{In fact, a given Dedekind zeta function only determines the product of the class number and the regulator.}. Using \sage \cite{sage}, we may compute a large amount of $a_n$ quickly. We introduce the vector:
\begin{equation}\label{eq.zetavector}
v_Z(F)=(a_1,\dots,a_{1000})\in\mathbb{Z}^{1000}.
\end{equation}

\begin{example}
For the dataset of $6206$ number fields with Galois group $C_8$ mentioned in Example~\ref{ex.6206}, the largest absolute value in the $1000 \times 6206$ $v_Z(F)$-entries $a_i$ is $109824$.
\end{example}

Given a finite set $\mathcal{F}$ of number fields $F$ and an invariant $c(F)$ for each $F\in\mathcal{F}$, we associate the following labeled dataset:
\begin{equation}\label{eq.zetadata}
\mathcal{D}_Z=\{v_Z(F)\rightarrow c(F):F\in\mathcal{F}\}.
\end{equation}
We may write $\zeta_F(s)$ as a product indexed by rational primes:
\begin{equation}
\zeta_F(s)=\prod_{p}E_p(s)^{-1} \ ,
\qquad
E_p(s):=\prod_{\mathfrak{p}|p}(1-N(\mathfrak{p})^{-s}) \ .
\end{equation}
If $\mathfrak{p}|p$ then $N(\mathfrak{p})$ has the form $p^a$ for $a\in\mathbb{Z}_{>0}$. Thus the product $E_p(s)$ is a polynomial in $p^{-s}$.

\begin{example}
Assume that $F$ is a quadratic extension of $\mathbb{Q}$. For a rational prime $p$, we have
\[
E_p(s)=
\begin{cases}
1-2p^{-2}+p^{-2s},&\text{ if } p\text{ is split,}\\
1-p^{-2s},&\text{ if } p\text{ is inert,}\\
1-p^{-s},&\text{ if } p\text{ is ramified.}
\end{cases}
\]
The splitting property of an unramified prime $p$ in $F$ is determined by the Legendre symbol  and the quadratic reciprocity law. See \cite[I, \S 8]{Neu} for more details. 
\end{example}

\subsection{Split primes}

For a number field $F$ that is Galois over $\mathbb Q$, write $\mathrm{Spl}(F)$ for the set of rational primes that split completely in $F$.  The Chebotarev density theorem implies that the set $\mathrm{Spl}(F)$ has density $1/[F:\mathbb Q]$, and it can be shown that
\[ 
F \subset K \qquad \Longleftrightarrow \qquad \mathrm{Spl}(F) \supset \mathrm{Spl}(K)
\]
for finite Galois extensions $F$ and $K$ over $\mathbb Q$ \cite[VII, \S 13]{Neu}. This shows that the set $\mathrm{Spl}(F)$ characterizes a Galois extension completely. 

For $i\in\mathbb{Z}_{\geq1}$, let $p_i$ denote the $i^{\mathrm{th}}$ rational prime. Given a number field $F$, we write
\begin{equation}\label{eq.deltai}
\delta_i=\begin{cases}
1,&\text{if }p_i\in\mathrm{Spl}(F),\\
0,&\text{otherwise.}
\end{cases}
\end{equation}
Except for finitely many primes, in order to calculate $\delta_i$ it suffices to reduce the defining polynomial $P(x)$ modulo $p_i$. If the reduction splits into a product of distinct linear factors then $\delta_i=1$; otherwise we have $\delta_i=0$. Using \sage \cite{sage}, it is possible to calculate a large number of $\delta_i$ quickly. Associated to $F$, we introduce the following binary vector:
\begin{equation}\label{eq.splittingvector}
v_B(F)=(\delta_1,\dots,\delta_{500})\in\{0,1\}^{500}.
\end{equation}
We note that the 500$^{\mathrm{th}}$ prime is $3571$. Given a finite set $\mathcal{F}$ of number fields and an invariant $c(F)$ for each $F\in\mathcal{F}$, we associate the following dataset:
\begin{equation}\label{eq.splittingdata}
\mathcal{D}_B=\{v_B(F)\rightarrow c(F):F\in\mathcal{F}\}.
\end{equation}

\subsection{Experimental strategy}\label{sec:strategy}

\begin{enumerate}

\item Let $\mathcal{F}$ be a finite set of number fields. The choice of $\mathcal{F}$ depends on the experiment. For example, $\mathcal{F}$ could be a random sample of degree $8$ extensions and discriminant less than some bound.

\item For a number field $F\in\mathcal{F}$, let $c(F)$ denote a certain invariant of interest. For example, $c(F)$ could be a binary digit (category) corresponding to whether or not $\mathrm{Gal}(F/\mathbb{Q})$ is abelian.

\item Generate datasets of the form $\mathcal{D}=\{v(F)\rightarrow c(F):F\in\mathcal{F}\}$, where $\mathcal{D}$ is as in \eqref{eq.polydata}, \eqref{eq.zetadata}, or \eqref{eq.splittingdata}.

\item Decompose $\mathcal{D}$ as a disjoint union $\mathcal{T} \sqcup \mathcal{V}$, where $\mathcal T$ is a training set and $\mathcal V$ is a validation set. We use various ratios for splits of $\mathcal T$ and $\mathcal V$ such as 80-20, 70-30 or 20-80 percentage-wise. As there is no significant difference in the results, we will not specify ratios for individual experiments.

\item Train a classifier on the set $\mathcal{T}$. In this paper we will use random forests and logistic regression, which we implement using \mth \cite{wolfram}.

\item For all unseen number fields $F\in\mathcal{V}$, ask the classifier to determine $c(F)$. We record the precision and confidence. Here, precision is defined to be the percentage agreement of the actual value with the one predicted by the classifier. As an extra check to minimize false positives and false negatives, the confidence in the form of Matthews' correlation coefficient \cite{matthews} is computed. Both precision and confidence are desired to be close to 1.

\end{enumerate}

\section{Degree, signature, and rank}\label{sec:signature}

We recall that the extension degree $[F:\mathbb{Q}]$ is equal to $n=r_1+2r_2$, and that  the rank of the unit group $\mathcal{O}_F^{\times}$ is $r=r_1+r_2-1$. The Dedekind zeta function vanishes to order $r$ at $s=0$. To perform the experiments in this section we downloaded datasets from \cite[Number~fields]{lmfdb}. The completeness of this data is documented at \cite[Completeness~of~number~field~data]{lmfdb}. 

\subsection{Experiment I: Extension degree}

Whilst the defining polynomial of a Galois extension clearly encodes the extension degree (as the degree of the polynomial), the same is not obviously true for zeta coefficients or split prime data. Datasets consisting of Galois extensions of $\mathbb{Q}$ with Galois group $C_4$, $C_6$ and $C_8$ are obtained from \cite[Number~Fields]{lmfdb} and thus a $3$-category label can be established in the form of \eqref{eq.splittingdata}:
\begin{equation}
\mathcal{D}_B=\{ (\delta_1,\dots,\delta_{500}) \rightarrow  c \} \ ,
\end{equation}
where $\delta_i \in \{0,1\}$ and $c=0,1,2$, say, according to which of the $3$ Galois groups the number field $F$ corresponds. As with all cases below, in order to balance the data, we sample around $6200$ in each category. We find that, when trained on split prime data, a {\it logistic regression classifier} is able to perform this $3$-way classification with precision $0.976$ and confidence $0.968$. Even better, when trained on zeta coefficient data, a {\it random forest classifier} performs the same classification with precision $0.999$ and confidence $0.998$. 

\subsection{Experiment II: Rank of unit group}

Recall that the signature $(r_1,r_2)$ determines the rank $r$ of $\mathcal{O}_F^{\times}$ through $r=r_1+r_2-1$, and $\mathrm{sgn}(\Delta_F)=(-1)^{r_2}$.

\begin{example}
If $[F:\mathbb{Q}]=2$, then the unit group $\mathcal{O}_F^{\times}$ has rank $1$ (resp. $0$) if $F$ is totally real (resp. imaginary), that is, that the signature of $F$ is $(2,0)$ (resp. $(0,1)$).  The number field $F$ has rank $1$ (resp. $0$) if and only if $\Delta_F>0$ (resp. $\Delta_F<0$). More generally, if $[F:\mathbb{Q}]$ has even degree, then the unit group has odd (resp. even) rank if and only if $\Delta_F>0$ (resp. $\Delta_F<0$). We note that if $[F:\mathbb{Q}]=6$ (resp. $8$), then the rank of $\mathcal{O}_F^{\times}$ is an integer in the set $\{2,3,4,5\}$ (resp. $\{3,4,5,6,7\}$).
\end{example}

We obtain, from \cite[Number~Fields]{lmfdb}, the datasets consisting of Galois extensions of $\mathbb{Q}$ with cyclic Galois group and signatures: $(2,0)$ and  $(0,1)$ for $C_2$, $(6,0)$ and $(0,3)$ for $C_6$, $(8,0)$ and $(0,4)$ for $C_8$. These signatures correspond to ranks $1,0,5,2,7,3$ respectively. Furthermore, we downloaded those with Galois group $D_4$ and signatures $(8,0)$ and $(0,4)$. This establishes datasets in the form of \eqref{eq.polydata} and \eqref{eq.zetadata}, {\it for each} choice of the 4 Galois groups, as
\begin{equation}
\begin{split}
\mathcal{D}_P &=\{ (c_0,\dots,c_{n-1}) \rightarrow r\} \ , \\
\mathcal{D}_Z &= \{(a_1,\dots,a_{1000}) \rightarrow r\} \ .
\end{split}
\end{equation}
Here, $c_i \in \IZ$ are the coefficients of the (monic) defining polynomial and $a_i$ are the first 1000 coefficients of the Dedekind zeta function. The rank $r$, conveniently, takes values with the binary categories for each Galois group, as indicated in Table \ref{t:rank}.

\begin{table}[h!]
\begin{center}
{\scriptsize \begin{tabular}{|c|c|c|c|c|}
\hline
Galois group & signature$(F)$ & rank$(\cO_F^{\times})$ &  $\mathcal D_P$ precision & $\mathcal D_P$ confidence\\ 
\hline
$C_2$ & \begin{tabular}{@{}l@{}}
                   (2,0) \\
                   (0,1) \\
                 \end{tabular} & 
\begin{tabular}{@{}l@{}}
                   1 \\
                   0 \\
                 \end{tabular} & $>0.99$ & $>0.99$ \\
\hline
$C_6$ & \begin{tabular}{@{}l@{}}
                   (6,0) \\
                   (0,3)  \\
                 \end{tabular} & 
\begin{tabular}{@{}l@{}}
                   5 \\
                   2 \\
                 \end{tabular} & 0.97 & 0.93\\
\hline
$C_8$ & \begin{tabular}{@{}l@{}}
                   (8,0) \\
                   (0,4)  \\
                 \end{tabular} & 
\begin{tabular}{@{}l@{}}
                   7 \\
                   3 \\
                 \end{tabular} & $>0.99$ & $>0.99$ \\
\hline
$D_4$ & \begin{tabular}{@{}l@{}}
                   (8,0) \\
                   (0,4) \\
                 \end{tabular} & 
\begin{tabular}{@{}l@{}}
                   7 \\
                   3 \\
                 \end{tabular}  & 0.98 & 0.95\\
\hline
\end{tabular} }
\end{center}
\caption{{\sf
Random forest classifier results on distinguishing ranks of the unit group for field extensions $F$ over $\IQ$ with fixed Galois groups.
}\label{t:rank}}
\end{table}

When trained on zeta coefficients using $\mathcal{D}_Z$ we found that all the standard classifiers, including neural-classifiers with convolutional networks, performed quite poorly.  In all the cases of Galois groups, the precision was around 0.6 or less.  It is interesting that this particular case requires so much effort without success whilst the majority of case which preform well amongst our experiments had good accuracies with many different classifiers. On the other hand, when trained on the defining polynomial coefficients using $\mathcal{D}_P$, the random forest classifier consistently performed the best and the precisions were $>0.99$, $0.97$, $>0.99$ and $0.98$ for the Galois groups as summarized in Table \ref{t:rank} with the corresponding confidences.

In fact, we can do more. Note that for our data, in each Galois group, the rank takes one of 2 possible values, which we can take to be 0 or 1 by appropriate labeling. This naturally makes one think of the logistic sigmoid function:
\begin{equation}
\sigma(z)=\frac{1}{1+\exp(-z)} \qquad \quad 
\raisebox{-30 pt}{\includegraphics[scale=0.2]{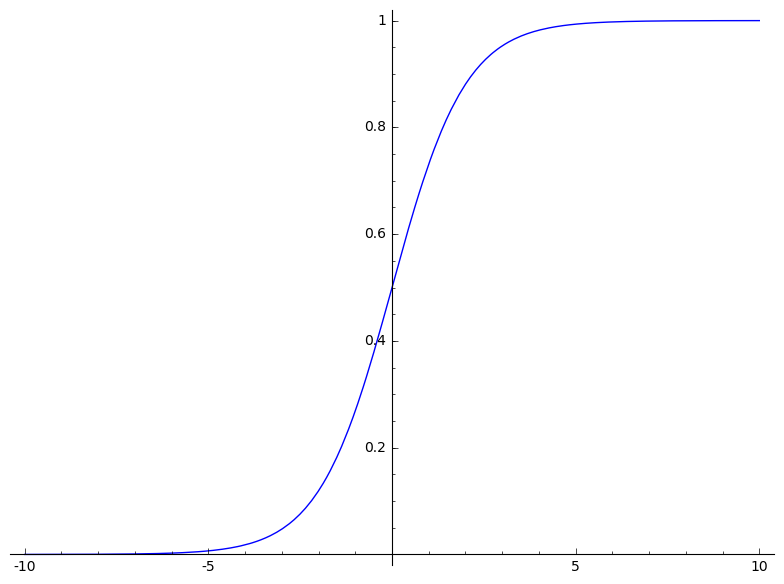}}
\end{equation}
which has range $[0,1]$ as shown in the graph above.

When trained on defining polynomial coefficients, we found that logistic regression performed well in predicting the rank of $\mathcal{O}_F^{\times}$, giving us an explicit and interpretable model. As with all regression, we are interested in finding a best fit to a function, here a sigmoid of the form:
\begin{equation}\label{eq.sigmoid}
\sigma(c_0 w_0+\cdots+c_{n-1} w_{n-1}+w_n), \ \ (w_0,\dots,w_n)\in\mathbb{R}^{n+1},
\end{equation}
where $c_i$, we recall,  are the coefficients of defining minimal polynomials for $F$ as in equation~\eqref{eq.minpoly}.
The parameters $w_i$ are to be optimized (fitted) by minimizing squared-mean-error. By rounding the above to the nearest integer, the function in equation~\eqref{eq.sigmoid} gives the value $0$ or $1$ corresponding to the two possibilities for the rank.
%
%
\begin{example}
In the case that $\mathrm{Gal}(F/\mathbb{Q})=C_6$,  regression by the logistic function yields around $94\%$ precision with best fit:
\begin{multline*}
c_0w_0+\cdots+c_{5}w_{5}+w_6 =-0.000169037 c_0 - 0.0000689721 c_1 - 0.000120625 c_2\\ 
- 0.00196535 c_3 - 0.058735 c_4 + 0.917924 c_5.
\end{multline*}
The accuracy of this model varies with ranks and $\Delta_F$. More precisely, the model predicts rank $2$ with accuracy $>0.91$ for almost all ranges of $|\Delta_F|$ occurring in the dataset with overall precision $0.987$. On the other hand, the model above performs poorly for rank $5$ fields with $|\Delta_F|<1.80 \times 10^9$
(around $77\%$ accuracy) and $|\Delta_F|>2.55 \times 10^{14}$ (around $60\%$ accuracy); the overall precision of rank 5 fields is $0.892$. 
\end{example}
We point out that what we did above should, strictly speaking, be called ``non-linear regression by the integer round of the logistic function''. The terminology {\it logistic regression}, though similar, is different in the probabilistic nature of the interpretation. Ordinarily, for discrete classification models such as our binary category problem, we fit the probability $p$ of the output being 0 or 1 to the logistic function.

\section{Galois groups of order 8}\label{sec:GalGroups}

Traditional algorithms for computing Galois groups are presented in \cite[Section~6.3]{Cohen}. In this section we will see that a classifier can distinguish between various Galois groups of order $8$, in a very efficient way. There are five possibilities for such groups, namely: $C_8$, $C_4\times C_2$, $C_2\times C_2\times C_2$, $D_4$, and $Q_8$ (the quaternion group). The last of these ($Q_8$) has about $1100$ occurrences on the LMFDB, so we do not use it for experimentation; the others range from around 6200 to 28000 cases. 
All of these can be obtained from  \cite[Number~Fields]{lmfdb}.

\subsection{Experiment III: Abelian vs. non-Abelian groups}\label{ss:abnonab}

Let us consider Galois extensions of $\mathbb{Q}$ with Galois group $C_8$ (resp. $D_4$). Note that $C_8$ is Abelian but $D_4$ is not. We establish a dataset of the form of \eqref{eq.polydata} as
\begin{equation}
\mathcal{D}_P = \{ (c_0,\dots,c_{7}) \rightarrow a \} \ ,
\end{equation}
where the input is the list of the 8 non-trivial coefficients $c_i$ of the minimal polynomial (the leading coefficient is always 1) and the output $a$ is 0 or 1 according to whether the Galois group is $C_8$ or $D_4$. We find that a random forest classifier was able to distinguish between these groups with precision $0.971$ and confidence $0.941$. Similarly, if we use the zeta coefficients $(a_1,\dots,a_{1000})$ as input along the lines of \eqref{eq.zetadata}, the random forest classifier achieves precision $0.973$ and confidence $0.947$.

On the other hand, when we use the split primes data $\mathcal D_B$ in \eqref{eq.splittingdata}, the random forest classifier yields precision $0.736$. 

\subsection{Experiment IV: Distinguishing between Abelian groups}

The above experiment showed that a classifier could distinguish between Abelian versus non-Abelian groups. We now ask: can a similar classifier perform the more refined distinction between different Abelian groups? Here, we have 3 abelain Galois groups of order $8$: $C_8$, $C_4\times C_2$, and $C_2\times C_2\times C_2$. Using the zeta coefficient data with the output being one of the 3 categories, we find that a random forest classifier was able to distinguish between these groups with precision $0.95472$ and confidence $0.932148$.

\section{Class numbers}\label{sec:classnumbers}

The Dedekind zeta function of an algebraic number field has a simple pole at $s=1$. At this pole, the residue is computed by the analytic class number formula which involves various arithmetic invariants including the class number\footnote{
Specifically, the class number formula dictates that the Dedekind zeta function has a simple pole at 1 and
\[
\lim\limits_{s \to 1} (s-1) \zeta_F(s) = \frac{2^{r_1} (2\pi)^{r_2} \mathrm{Reg}_F h_F}{w_F \sqrt{|\Delta_F|}} \ ,
\]
where, in addition to the nomenclature in \S\ref{s:nomen}, $\mathrm{Reg}_F$ is the regulator, and $w_F$ is the number of roots of unity in $F$.}. 
Algorithms for computing the class numbers of general number fields are given in \cite[Section~6.5]{Cohen}. For the special case of quadratic extensions, see \cite[Sections~5.2,~5.6]{Cohen}. 

\subsection{Experiment V: Real quadratic fields}

In the discriminant range $0<|\Delta_F|\leq10^6$, one finds $83464$ real quadratic number fields with class number $1$, and $83324$ real quadratic number fields with class number $2$ \cite[Number~Fields]{lmfdb}, and the list is complete for this discriminant range.

Whilst there are many thousands of examples of real quadratic number fields with larger class number, the sample size varies from case to case. In order to avoid biases in our datasets, we simply focus on the binary classification problem of distinguishing real quadratic fields with class number $1$ from those with class number $2$. 
Thus, we have datasets, using \eqref{eq.polydata}, \eqref{eq.zetadata} and \eqref{eq.splittingdata},
\[
\mathcal{D}_P =\{ (c_0,c_1) \rightarrow c \},\quad
\mathcal{D}_Z  =\{ (a_1,\dots,a_{1000}) \rightarrow c \}, \quad 
\mathcal{D}_B  =\{ (\delta_1,\dots,\delta_{500}) \rightarrow c \},
\]
where $c_i$ are the coefficients of the minimal polynomial, $a_i$ are the first 1000 zeta coefficients, and $\delta_i$ are 0 or 1 according to whether the $i$th rational prime splits completely or not.
Here $c = 1$ or $2$ is the class number.
Note that we have a fairly balanced dataset with around $80,000$ each of class number 1 and 2.
When trained on zeta coefficient data, the random forest classifier yielded the best precision of $0.96$ with confidence $0.92$. This experiment is summarized in the table below. On the other hand, when trained on defining polynomial data or split primes data, no standard classifier was able to distinguish between class numbers $1$ and $2$ with precision greater than $0.60$.

Next, we try something more drastic.
Consider the discriminant range $10^6<\Delta_F<2\times10^6$, in which we find $75202$ real quadratic number fields with class number $1$ and $80217$ with class number $2$. According to \cite[Number~Fields]{lmfdb}, the list is complete for this discriminant range.

Can a classifier be trained within $|\Delta_F|$ of a certain range and {\em extrapolate} to a larger $|\Delta_F|$ range?
If so, this would strengthen even further our notion that machine-learning has found some underlying pattern. 
We applied the classifier trained on the previous datasets, i.e., real quadratic fields with discriminant less than one million, to this new discriminant range. The result was precision $0.92$ with confidence $0.86$. It seems that the classifier is able to {extrapolate} from data of smaller discriminant.
We tried the same for discriminants between $2$ million and $3$ million, again with the same classifier trained on the data of discriminants smaller than one million. There are $18383$ with class number $1$, and $19827$ with class number $2$. The result was precision $0.91$ with confidence $0.84$. 
The results are summarized in Table \ref{t:deltaRange}.

\begin{table}
\begin{center}
{\scriptsize \begin{tabular}{|c|c|c|c|c|}
\hline
Discriminant range&$h_F$&$\#\{F\}$&Precision&Confidence\\
\hline
$[1,1\times10^6]$& \begin{tabular}{@{}c@{}}
                     1\\ 2\\ 
                 \end{tabular} &\begin{tabular}{@{}c@{}}
                     83464\\ 83324 \\
                 \end{tabular} &0.96&0.92\\
\hline
$[1\times10^6,2\times10^6]$& \begin{tabular}{@{}c@{}}
                     1\\ 2\\ 
                 \end{tabular} &\begin{tabular}{@{}c@{}}75202\\80217 \\
\end{tabular}&0.92&0.86\\
\hline
$[2\times10^6,3\times10^6]$& \begin{tabular}{@{}c@{}}
                     1\\ 2\\ 
                 \end{tabular} &\begin{tabular}{@{}c@{}}18383\\ 19827\\
\end{tabular}&0.91&0.84\\
\hline
\end{tabular}
}
\end{center}
\caption{{\sf
A summary of the precision and confidence of the random forest classifier trained on zeta coefficients of real quadratic fields with discriminant between one and one million with class number 1 or 2. The number $\#\{F\}$ is the cardinality of the set containing real quadratic fields with discriminant and class number as specified. 
}}
\label{t:deltaRange}
\end{table}

\begin{remark}
It is known that there is a finite set of imaginary quadratic fields with class number $1$, viz., this is the list of
$\IQ[\sqrt{d}]$ for the Heegner numbers $d= -1, -2, -3, -7, -11, -19, -43, -67, -163$.
This set has far too few examples for a machine to learn.  Motivated by \cite{CL}, we tested to see if any classifier could distinguish between class numbers divisible by $3$ and class numbers not divisible by $3$: The methods of zeta coefficients and polynomial coefficients gave both precision around $0.51$, which is as good as randomly guessing.
As always, divisibility and other patterns in primes seem very difficult to be machine-learned (cf. \cite{He:2018jtw}).
\end{remark}

\subsection{Experiment VI: Quartic and sextic fields}

We say a degree $4$ number field is {\em bi-quadratic} if it has Galois group $C_2\times C_2$. From \cite[Number~Fields]{lmfdb} we downloaded the dataset of bi-quadratic Galois extensions of $\mathbb{Q}$ with class number $1$ and $2$. To get a balanced dataset, we randomly chose $6100$ number fields for each class number. When trained on zeta coefficient data, we found that a logistic regression classifier could distinguish class number 1 from class number 2 with precision $0.819$ with confidence $0.640$. We suspect that the performance could have been better with a larger set of data. 

In degree $6$, the generic Galois group is $S_3$. From \cite[Number~Fields]{lmfdb}, we downloaded the dataset consisting of degree $6$ Galois extensions of $\mathbb{Q}$ with Galois group $S_3$ and class number in the set $\{1,2,3,4,6,8,9\}$. The class numbers $5$ and $7$ were excluded on the grounds that there are too few data points on the LMFDB (less than $400$ each), whereas the others have at least $1150$ points; precise counts are given in Table \ref{t:46}.

\begin{table}
\begin{center}
{\scriptsize
\begin{tabular}{|c||c|c|c|c|c|c|c|}
\hline
$h_F$&1&2&3&4&6&8&9\\
\hline
$\#\{F\}$& 7436 & 8680 & 1917 & 8165 & 1158  & 4230 & 2700 \\
\hline
\end{tabular}
}
\end{center}
\caption{{\sf Frequency of $S_3$-extensions $F$ of $\mathbb{Q}$ with class number $h_F$.}}
\label{t:46}
\end{table}

When trained on zeta coefficient data consisting of randomly chosen $1150$ data points from each class number, we found no machine-learning approach was able to perform the corresponding $7$-way classification with precision more than $0.38$. Furthermore, when trained on a dataset consisting of $7400$ data points from each of the  class numbers $1,2,4$, the best precision given by a random forest classifier was $0.605$ with confidence $0.415$.   
Even with polynomial coefficient data, we found that no classifier could distinguish between class number $1$ and class number $2$ with precision greater than $0.64$. 
%
This poor performance might well be due to the small size of the dataset
(recall that the size of the dataset for real quadratic fields of class number 1 is larger than 83000).

\section{Outlook}

We conclude with a brief discussion of future experimental and mathematical projects.

As mentioned in the introduction, it is unknown whether or not there are infinitely many real quadratic fields of class number $1$. It would be very interesting to investigate how a machine is able to distinguish such fields. If the criteria under which the classifier predicts class number $1$ are satisfied infinitely often, then there could be scope for developing a new heuristic for or a new approach to this open problem. Furthermore, we note that the machine continues to make accurate predictions for real quadratic fields with discriminant outside the range of the training data. Perhaps this extrapolation offers a clue towards future progress.

The class number is subject to the analytic class number formula, which computes the residue of the Dedekind zeta function at its pole. The analytic class number formula can be compared to the famous BSD conjecture in the sense that both concern the leading terms of zeta functions at special points. In a forthcoming paper, we will examine whether or not a machine can be trained to predict the vanishing orders of elliptic $L$-functions and other invariants appearing in their Taylor expansions.


{\small 
Yang-Hui He {\sf hey@maths.ox.ac.uk} \\
Department of Mathematics, City, University of London, EC1V 0HB, UK;\\
Merton College, University of Oxford, OX14JD, UK;\\
School of Physics, NanKai University, Tianjin, 300071, P.R.~China

Kyu-Hwan Lee {\sf khlee@math.uconn.edu} \\
Department of Mathematics, University of Connecticut, Storrs, CT, 06269-1009, USA

Thomas Oliver {\sf Thomas.Oliver@nottingham.ac.uk} \\
School of Mathematical Sciences, University of Nottingham, University Park, \\
Nottingham, NG7 2QL, UK
}

\end{document}